\documentclass{amsart}
\usepackage{palatino}
\usepackage{amsthm, amsmath}
\usepackage{amssymb} 
\usepackage{amscd}
\usepackage[all]{xy}
\usepackage[breaklinks,bookmarksopen,bookmarksnumbered]{hyperref}
\usepackage{geometry}
\theoremstyle{plain}
\newtheorem*{thm}{Theorem}
\newtheorem{lem}{Lemma}
\newtheorem{prop}{Proposition}
\newtheorem*{prop*}{Proposition}
\newtheorem{cor}{Corollary}
\theoremstyle{remark}

\newtheorem{ex}{Examples}

\def\pr{\noindent\textit{Proof} : }

\def\rond{\kern 1pt{\scriptstyle\circ}\kern 1pt}

\newcommand\Sp{\operatorname{Sp}}
\newcommand\PSp{\operatorname{PSp}}
\newcommand\GL{\operatorname{GL}}
\newcommand\SL{\operatorname{SL}}
\newcommand\PSL{\operatorname{PSL}}
\newcommand\PGL{\operatorname{PGL}}
\newcommand\ed{\operatorname{ed}}
\newcommand\rd{\operatorname{rd}}

\newcommand\SO{\operatorname{SO}}

\def\Z{\mathbb{Z}}
\def\Q{\mathbb{Q}}

\def\C{\mathbb{C}}
\def\P{\mathbb{P}}

\def\F{\mathbb{F}}

\def\A{\mathfrak{A}}

\def\iso{\vbox{\hbox to .8cm{\hfill{$\scriptstyle\sim$}\hfill}
\nointerlineskip\hbox to .8cm{{\hfill$\longrightarrow $\hfill}} }}

\begin{document}
\title[Finite simple groups of small essential dimension]{Finite simple groups of small essential dimension}
\author[Arnaud Beauville]{Arnaud Beauville}
\address{Laboratoire J.-A. Dieudonn\'e\\
UMR 7351 du CNRS\\
Universit\'e de Nice\\
Parc Valrose\\
F-06108 Nice cedex 2, France}
\email{arnaud.beauville@unice.fr}
 
\date{\today}
 
\begin{abstract}
We discuss the notion of essential dimension of a finite group (over $\C$) and explain its relation with birational algebraic geometry. We show how this leads to a (partial) classification of  simple finite groups of essential dimension $\leq 3$.
\end{abstract}

\maketitle 

\section*{Introduction}

The \emph{essential dimension} of a finite group $G$ has been introduced in the seminal paper \cite{BR}. We have to refer to that paper for  motivation (see also the survey paper \cite{Bgaz}); in a very informal way, the essential dimension of $G$ (over a field $k$) is the minimum number of parameters needed to define a general Galois extension $L/K$ with Galois group $G$ and $K\supset k$. 

Since then the notion of essential dimension has been put in a much larger context (see \cite{Merk} for a recent survey). In this note we want to go back to the original problem, in the simplest case where $k$ is the field of complex numbers. In that case the definition is quite concrete, and the computation of the essential dimension becomes an interesting problem  of classical algebraic geometry.  We will explain how one can use 
 results on the birational geometry of unirational varieties  to classify the simple finite groups of essential dimension $\leq 2$, and obtain partial results in dimension 3.

\section{Definitions and basic properties}

Let $G$ be a finite group, and $X$ a complex
algebraic variety  with a faithful action of $G$. We will say that $X$ is \emph{linearizable} if there exists a complex
representation $V$ of $G$ and a rational dominant $G$-equivariant map $V\dasharrow X$. %(such a map is called a {\it compression} of $V$). 
The \emph{essential dimension} $\ed(G)$ of $G$ (over $\C$) is the
 minimal dimension of all linearizable $G$-varieties.  
 
 \smallskip	
 Let us mention immediately three obvious consequences of the definition, which we will use frequently in the sequel :
 
 $\bullet\ $ One has $\ed(G)\geq 0$, and $\ed(G)=0$ if and only if $G=\{1\} $;
 
 $\bullet\ $ If $H$ is a subgroup of $G$, $\ed(H)\leq \ed(G)$.
 
 $\bullet\ $ Let $u:X\dasharrow Y$ be a $G$-equivariant birational map; then $X$ is  $G$-linearizable if and only if $Y$ is. In particular, to compute $\ed(G)$ it suffices to consider smooth projective varieties.	
\begin{ex}\label{ex}
$a)$ A vector space $V$ with a faithful linear action of $G$ is obviously $G$-linea\-rizable, hence $\ed(G)\leq \dim (V)$. Therefore, if we denote by $\rd(G)$ the minimal dimension  of a faithful representation of $G$ (``\emph{representation dimension}'' of $G$), we have $\ed(G)\leq \rd(G)$.

$b)$ In the situation of $a)$, the group $G$ acts on the projective space $\P(V)$, and the rational map $V\dasharrow \P(V)$ is $G$-equivariant. The action of $G$ on $\P(V)$ is faithful if and only if $G\subset \GL(V)$ does not contain any nontrivial homothety; this is the case in particular if the center of $G$ is trivial. In this situation we have $\ed(G)\leq \rd(G)-1$.

$c)$ Here is a more elaborate example. The permutation action of the symmetric group $\mathfrak{S}_n$ on $\C^n$  extends to $(\P^1)^n$. On the other hand the group $H=\PGL(2,\C)$ acts also on $(\P^1)^n$; the Zariski open subset $(\P^1)^n_{\neq }$ of points with distinct coordinates is stable under both actions. For $n\geq 3$ $H$ acts freely on $(\P^1)^n_{\neq }$; the quotient  $X_n=(\P^1)^n_{\neq }/H$ is an algebraic variety of dimension $n-3$. The action of $\mathfrak{S}_n$ on  $(\P^1)^n_{\neq }$ commutes with that of $H$, hence induces an action on $X_n$, which is faithful for $n\geq 5$. We conclude
that \emph{the essential dimension of $\mathfrak{S}_n$ is $\leq n-3$} for $n\geq 5$.
 \end{ex}

The following result will provide our main tool to prove that a variety is \emph{not} $G$-linearizable :
 \begin{prop}[\cite{KS}]\label{ks}
Let $A$ be a finite abelian group, $X$ an $A$-linearizable  projective variety. There is a point of $X$ fixed by $A$. 
\end{prop}
The proof is so simple that we cannot resist to to copy it from \cite{KS}. Since a linear action fixes the origin, the Proposition follows from a more general result :
\begin{lem}
Let  $X$, $Y$ be two $A$-varieties, with $X$ smooth and $Y$ proper, $f:X\dasharrow Y$  a rational $A$-equivariant map. If $X$ has a point fixed by $A$, so does $Y$.
\end{lem}
\pr The proof is by induction on $\dim(X)$, the case $\dim(X) = 0$ being clear. 

Let $x\in X$ be a  point fixed by $A$; let $B_x(X)$ be the blow up of $x$ in $X$, and $E$ the exceptional divisor. The action of $A$ extends to $B_x(X)$ and $E$. Since $Y$ is proper, the rational map $B_x(X)\dasharrow Y$ is defined outside a subset of codimension $\geq 2$, so it induces a $A$-equivariant rational map $E\dasharrow Y$. Since an abelian group acting on a projective space has always a fixed point, $A$  fixes a point of $E$, hence of $Y$ by the induction hypothesis.\qed

\begin{cor}\label{eda}
We have $\ed(A)=\rd(A)$. In particular, the essential dimension of $(\Z/n)^r$ est $r$.
\end{cor}
\pr Let $X$ be a smooth $A$-linearizable projective variety, $x$ a point of $X$ fixed by $A$. The group $A$ acts on the tangent space $T_x(X)$, and this action is  isomorphic in a neighborhood of $x$ to the action of $A$ on $X$. Therefore the representation of $A$ on $T_x(X)$ is faithful; thus $\rd(A)\leq \dim(X)$, hence $\rd(A)\leq \ed(A)$. The opposite inequality is obvious (example --). Finally the equality $\rd((\Z/n)^r)=r$ is  an easy exercise.\qed
 
 \medskip	
 
The equality $\ed(G)=\rd(G)$  holds more generally for a $p$-group \cite{KM}; the proof uses much more sophisticated techniques.
 
 \smallskip	
Let us conclude this section with the list of groups of essential dimension $1$ :
\begin{prop}[\cite{BR}] 
The finite groups of essential dimension $1$ are the cyclic groups and the dihedral groups $D_n$ for $n$ odd.\label{ed1}
\end{prop}
\pr We have $\rd(\Z/m)=1$, and $\rd(D_n)=2$; when $n$ is odd the center of $D_n$ is trivial. Thus $\ed(\Z/m)=\ed(D_n)=1$ by Examples \ref{ex}.$a)$ and $b)$.

Let $G$ be a finite group with $\ed(G)=1$, and $X$ a $G$-linearizable smooth projective curve. Then $X\cong \P^1$, so $G$ is a subgroup of $\PGL_2(\C)$, hence isomorphic to $\Z/n,D_n,\mathfrak{A}_4,\mathfrak{S}_4$ or $\mathfrak{A}_5$. Now except $\Z/m$ and $D_n$ for $n$ odd, all these groups contain a copy of $(\Z/2)^2$. We conclude with 
Corollary \ref{eda}.

\medskip	
\section{Groups of essential dimension 2}

The groups of essential dimension 2 have been classified in \cite{D2}; the list is already quite large. We will restrict ourselves to the class of simple groups.

\begin{prop}
The simple finite groups of essential dimension $2$ are $\mathfrak{A}_5$ and $\PSL_2(\F_7)$. 
\end{prop}
Of course this follows immediately from \cite{D2}; but the proof is significantly easier for simple groups. 

\smallskip	
\pr We have $\ed(\mathfrak{A_5})= \ed(\mathfrak{S_5})=2$ by Example \ref{ex}.$c)$ and  Proposition \ref{ed1}; the group $\PSL_2(\F_7)$ has a faithful representation of dimension 3 (which can be realized as $H^0(C,K_C)$, where $C$ is the Klein quartic), hence it has essential dimension  2 by Example \ref{ex}.$b)$ and Proposition \ref{ed1}. 

 Let $G$ be a finite group with $\ed(G)=2$, and $X$ a $G$-linearizable smooth projective 
surface. By Castelnuovo's theorem $X$ is rational, so $G$ is a finite subgroup of the Cremona group $\mathrm{Cr}_2=\mathrm{Bir}(\P^2)$. The classification of these subgroups has been worked out in the 19th century (Kantor, Wiman), and completed in \cite{DI}. Let us state the result in the case of interest for us, namely that of simple groups; here again, the proof  is much easier in that case than for general groups.

\begin{thm}[\cite{DI}]
The simple  finite subgroups of $\mathrm{Cr}_2$ are cyclic or isomorphic to $\mathfrak{A}_5$, $\mathfrak{A}_6$ or $\PSL_2(\F_7)$.\qed
\end{thm}

Thus our task now is to eliminate the group $\mathfrak{A}_6$. This group appears only once in the list, as a group of automorphisms of $\P^2$ (the so-called Valentiner group).
The inverse image $\tilde{\mathfrak{A}}_6 $ of $\mathfrak{A}_6$ in $\SL(3)$ is a central extension of $\mathfrak{A}_6$ by $\Z/3$. 

One way of describing this extension is to view $\mathfrak{A}_6$ as the subgroup of $\PGL(3,\F_4) $ preserving the  set formed by the 6 points
\[(1,0,0)\ ,\ (0,1,0)\ ,\ (0,0,1)\ ,\ (1,1,1)\ ,\ (1,\alpha ,\beta )\ ,\ (1,\beta ,\alpha )\ ,\]
where $\F_4=\{0,1,\alpha ,\beta \}$;
then  
$\tilde{\mathfrak{A}}_6 $ is the pull back of ${\mathfrak{A}}_6 $ in $\GL_3(\F_4)$. 
The elements 
\[ u = \begin{pmatrix}
1 & 0 & 0\\ 0 & \alpha & 0 \\ 0 & 0 & \beta 
\end{pmatrix}\qquad 
v = \begin{pmatrix}
0 & 1 & 0\\ 0 & 0 & 1 \\ 1 & 0 & 0 
\end{pmatrix}
\]of $\GL_3(\F_4)$ belong to $\tilde{\mathfrak{A}}_6 $, and their commutator is a nontrivial element of its center. We conclude with 
 the following lemma :
\begin{lem}\label{commut}
Let $V$ be a complex vector space, 
$G$  a finite subgroup of $\PGL(V)$, $\tilde{G} $ its inverse image in $GL(V)$. Assume that there exist elements $u,v$ of $\tilde{G} $ such that their commutator $(u,v)$ is a nontrivial homothety. 
Then $\P(V)$ is not $G$-linearizable. 
\end{lem}
\pr  There is no line in $V$ stable under $u$ and $v$, since otherwise $u$ and $v$ would commute on that line. Therefore the images of $u$ and $v$ in $\PGL(V)$  do not fix a common point of $\P(V)$. Since they commute, Proposition \ref{ks} shows that $\P(V)$ is not $G$-linearizable.\qed

\smallskip	
\noindent\emph{Remark}$.-$ Assume moreover that $G$ is simple, and that the action of $G$ on $\P(V)$ does not come from a linear action on $V$. Then it follows from \cite{Bl} that the 
hypothesis of the lemma is always satisfied. Thus $\P(V)$ is never $G$-linearizable in that case.

\medskip	
\section{Groups of essential dimension 3}
\subsection{Prokhorov's list}
The classification of all finite groups of essential dimension 3 is definitely out of reach at this moment. On one hand the classification of finite subgroups of the Cremona group $\mathrm{Cr}_3$ seems untractable; moreover, a $G$-linearizable threefold is unirational, but there is no reason to expect it to be rational.

However when we restrict our attention to simple groups, using a remarkable theorem of Prokhorov we get a partial result :
\begin{prop}\label{ed3}
The simple groups of essential dimension $3$ are $\mathfrak{A}_6$ and possibly $\PSL_2(\F_{11})$.
\end{prop}

\pr We have indeed $\ed(\mathfrak{A}_6)=3$ by Example \ref{ex}.$c)$ and the previous classification. The group $\PSL_2(\F_{11})$ admits a faithful (= nontrivial) representation of dimension 5, so  its essential dimension is 3 or 4 by Example \ref{ex}.$b)$; the exact value is not known (see (\ref{fin})).

 Let $G$ be a finite simple group with $\ed(G)= 3$. By definition there exists a linearizable
projective $G$-threefold $X$. This implies in particular that $X$ is rationally connected. Such pairs $(G,X)$ have been classified in \cite{P}: up to birational equivalence, we have the following possibilities for $G$ :
\begin{enumerate}
\item $G=\mathfrak{A}_5$, $\mathfrak{A}_6$,  $\PSL_2(\F_{7})$, $\PSL_2(\F_{11})$;
\item $G=\mathfrak{A}_7$, $\PSp_4(\F_3)$, $\SL_2(\F_8)$.
\end{enumerate}
\par The groups of the first row have already been dealt with. We will deal separately with the 3 remaining cases.

\subsection{The group $\mathfrak{A}_7$}
\begin{prop}[\cite{D1}]
The essential dimension of $\mathfrak{A}_7$ is $4$.
\end{prop}
\pr The group $\mathfrak{A}_7$ appears twice in Prokhorov's list :

$\bullet\ $ $\A_7$ acts by permutation of coordinates on the variety $X$ given by $\sum X_i=\sum X_i^2=\allowbreak \sum X_i^3=0$ in $\P^6$;

$\bullet\ $ $\A_7$ embeds into $\mathrm{PGL}(4,\C)$, hence acts on  $\P^3$.

\smallskip
In the first case, one checks easily that the subgroup $(\Z/2)^2\times \Z/3\subset \A_4\times \A_3\subset \A_7$  has no fixed point on $X$, so that $X$ is not $\A_7$-linearizable by Proposition \ref{ks}. 

In the second case, let us consider the double coverings of complex Lie groups
\[\mathrm{SL}(4)\longrightarrow \mathrm{SO}(6)\longrightarrow \mathrm{PGL}(4)\]
deduced from the isomorphism $\C^6\cong \wedge^2\C^4$. The standard  representation   $\A_7\subset\mathrm{SO}(6)$, composed with the second arrow, gives the embedding of $\A_7$  into $\mathrm{PGL}(4)$. The inverse image of
 $\A_7$ in $\mathrm{SL}(4)$  appears as a central extension 
\[ 1\rightarrow \{\pm I\} \longrightarrow  \tilde{\mathfrak{A}} _7 \longrightarrow  \A_7\rightarrow 1\ .\]
To apply Lemma \ref{commut}, we need to show that the central element $-I$ of $\tilde{\mathfrak{A}} _7$ is a commutator. Consider the subgroup
 $\A_4$ of $\A_7$ fixing the last 3 letters. We have a commutative diagram
\[ \xymatrix{\tilde{\mathfrak{A}} _4 \ar@{^{(}->}[rr]\ar[dd]\ar[rd] && \tilde{\mathfrak{A}} _7\ar[dr]\ar[dd] &\\
& \mathrm{SL}(2)\ar[dd]\ar@{^{(}->}[rr]&& \mathrm{SL}(4)\ar[dd] \\
\A_4 \ar@{^{(}->}[rr]\ar[rd] && \A_7\ar[dr] &\\
& \mathrm{SO}(3)\ar@{^{(}->}[rr]&& \mathrm{SO}(6)\\
}\] so it suffices to check that $-I$ is a commutator in $\tilde{\mathfrak{A}} _4$. Since 
  $\mathrm{SL}(2)$   contains no element of order 2 except $-I$, the Klein subgroup
   $(\Z/2)^2\subset\A_4$ lifts to the quaternion group $\Q_8=\{\pm 1,\pm i,\pm j,\pm k\} $, and we have $(i,j)=-1$. Thus $-I$ is a commutator in $\tilde{\mathfrak{A}} _7$, and Lemma \ref{commut} shows that $\P^3$ is not $\mathfrak{A}_7$-linearizable.
  \qed

\smallskip	
\subsection{The group $\PSp (4,\F_3)$}
\begin{prop}
The essential dimension of $\PSp (4,\F_3)$ is $4$.
\end{prop}
\pr  
The group $\Sp (4,\F_3)$ has a linear representation on the space $W$ of functions on $\F_3^2$, the \emph{Weil representation}, for which we refer to \cite{AR}, Appendix I. This representation splits as $W=W^+\oplus W^-$, the spaces of even and odd functions; we have $\dim W^+=5$, $\dim W^-=4$. The central element $(-I)$ of $\Sp (4,\F_3)$ acts on $W$ by ${}^{(-I)\!}F\,(x)= F(-x)$, hence it acts trivially on $W^+$, and as $-\mathrm{Id}$ on $W^-$. Thus we get 
 a faithful representation of $\PSp (4,\F_3)$ on $W^+$,
 with a compression to $\P(W^+)\cong\P^4$, hence $\ed (\PSp (4,\F_3))\leq 4$.
\par To prove that we have equality, we observe\footnote{I am indebted to A. Duncan for this observation.} that $\PSp (4,\F_3)$ contains a subgroup isomorphic to $(\Z/2)^4$. One way to see this is to use the classical isomorphism of $\PSp (4,\F_3)$ with $\SO^+(5,\F_3)$, the kernel of the spinor norm $\SO(5,\F_3)\rightarrow \F_3^*\cong\{\pm 1\}$. Essentially by definition (see \cite{Bo}, \S 9, no. 5), this group contains  the transformations $\sigma _v:x\mapsto -x+2(x.v)v\ $ for each length 1 vector $v$; when $v$ runs over the elements of an orthonormal basis, the $\sigma _v$ span a subgroup of $\SO^+(5,\F_3)$ isomorphic to $(\Z/2)^4$. By \cite{BR} we have \[ \ed(\PSp (4,\F_3))\geq \ed((\Z/2)^4)=4\ . \eqno{\qed}\]

 \medskip

\def\go{\mathbb{G}_{\mathrm{iso}}(4,U)}

\subsection{The group $\mathrm{SL}_2(\F_8)$}
\begin{prop}
The essential dimension of $\SL_2(\F_8)$ is $\geq 4$.
\end{prop}
The group $\SL_2(\F_8)$ has a representation of dimension 7, hence its essential dimension is $\leq 6\ $ -- we do not know its precise value.

\medskip
\pr The group $\SL_2(\F_8)$ acts on a rational Fano threefold $X\subset \P^8$ in the following way \cite{P}. Let $U$ be an irreducible 9-dimensional representation of $\SL_2(\F_8)$; there exists a non-degenerate invariant quadratic form $q$ on $U$, unique up to a scalar. Then  $SL_2(\F_8)$ acts on
the orthogonal Grassmannian $\go$ of 4-dimensional isotropic subspaces of $U$. This Grassmannian admits a $O(q)$-equivariant embedding into $ \P^{15}$, given by the half-spinor representation \cite{Mu}. The threefold $X$ is the intersection of $\go$ with a subspace $\P^8\subset \P^{15}$ invariant under $\SL_2(\F_8)$. 

\par Let $N\subset \SL_2(\F_8)$ be the subgroup of matrices $\begin{pmatrix}
1 & a \\
0  &  1\\
\end{pmatrix}$, $a\in\F_8$.
We will show that $N$ has no fixed point in $\go$, and therefore in $X$. 

\par Let $\chi _U$ be the character of the representation $U$. We have $\chi _U(n)=1$ for $n\in N$, $n\neq 1$ (see for instance \cite{C}, 2.7). It follows that $U$ restricted to $N$ is the sum of the regular representation and the trivial  one; in other words, as a $N$-module we have
\[U=\C_1^2\oplus \sum_{\lambda \in\hat{N}\smallsetminus\{1\}}  \C_\lambda \ ,\]
where $\C_\lambda $ is the one-dimensional representation associated to the character $\lambda $. The subspaces $\C_\lambda  $ and $\C_\mu  $ for $\lambda \neq \mu $ are orthogonal for $q$; 
 since $q$ is non-degenerate, its restriction to each $\C_\lambda $ ($\lambda \neq 1$) and to  $\C_1^2$ must be non-degenerate.

\par Now any vector subspace $L\subset U$ stable under $N$ must be the sum of some of the $\C_\lambda $, for $\lambda \neq 1$, and of some subspace of $\C_1^2$; this implies that $L$ cannot be isotropic as soon as $\dim L\geq 2$. Hence $N$ has no fixed point on $\go$, and $X$ is not linearizable by Proposition \ref{ks}.\qed

\medskip
This finishes the proof of Proposition \ref{ed3}.\qed

\smallskip
\subsection{About $\mathrm{PSL}_2(\F_{11})$}\label{fin}

 According to \cite{P} there are two rationally connected threefolds with an action of $\PSL_2(\F_{11})$, the Klein cubic $X^\mathrm{k}\subset \P^4$ given by $\sum_{i\in \Z/5}X_i^2X_{i+1}=0$ and a Fano threefold  $X^\mathrm{a}\subset\P^9$ of degree 14, birational to $X^\mathrm{k}$. The group $\PSL_2(\F_{11})$ has order $660=2^2 .3.5.11$; its abelian subgroups are  cyclic, except the  2-Sylow subgroups which are isomorphic to $(\Z/2)^2$. A finite order automorphism of a rationally connected variety has always a fixed point (for instance by the holomorphic Lefschetz formula); one checks easily that a 2-Sylow subgroup of $\PSL_2(\F_{11})$ has a fixed point on both $X^\mathrm{k}$ and $X^\mathrm{a}$. So Proposition \ref{ks} does not apply, and another approach is needed.
In [D-R] the authors show that the equality $\ed(\mathrm{PSL}_2(\F_{11}))=3$ would follow from a  conjecture of Cassels and Swinnerton-Dyer on the existence of rational points on cubic hypersurfaces.

\end{document}